\documentclass[12pt]{article}
\usepackage{bbm}
\usepackage{amsfonts}
\usepackage{mathrsfs}

\usepackage{amssymb}
\usepackage{amsmath}

\renewcommand{\phi}{\varphi}

   \newcommand{\PX}{{\mathbb P}}

\newcommand{\cF}{{\cal F}}

\newtheorem{theorem}{Theorem}
\newtheorem{lemma}{Lemma}

\newtheorem{example}{Example}
\newtheorem{definition}{Definition}

\newtheorem{remark}{Remark}
\begin{document}

\title{Random Chain Recurrent Sets\\ for Random Dynamical Systems
\footnote{This research was partly supported by the NSF grant
0620539, the Cheung Kong Scholars Program and the K. C. Wong
Education Foundation.}}

   \author{Xiaopeng Chen$^1$ and Jinqiao Duan$^{2}$ \\
1. School of Mathematics and Statistics\\
Huazhong University of Science and Technology \\
  Wuhan 430074, China  \\
 \emph{E-mail: chenxiao002214336@yahoo.cn} \\
2. Department of Applied Mathematics\\ Illinois Institute of Technology \\
  Chicago, IL 60616, USA \\ \emph{E-mail: duan@iit.edu} \\
 }

\date{\today}
\maketitle

\emph{ \bf Abstract.} It is known by the Conley's theorem  that the
chain recurrent set $CR(\varphi)$ of a
 deterministic flow $\varphi$ on a compact metric space  is the
complement of the union of sets $B(A)-A$, where $A$ varies over
the collection of attractors and $B(A)$ is the basin of attraction
of $A$.  It has recently been shown that a similar decomposition
result holds for random dynamical systems on noncompact separable
complete metric spaces, but under a so-called \emph{absorbing
condition}. In the present paper, the authors introduce a notion
of random chain recurrent sets for random dynamical systems, and
then prove   the random   Conley's theorem on noncompact separable
complete metric spaces \emph{without} the absorbing condition.

\textbf {Keywords}: Chain recurrent sets; attractors; Conley's
theorem; random dynamical systems; cocycle

 \textbf{Mathematics Subject Classification}: 37H99,
 37B20, 37B25,37B35

\section{Introduction}

 Chain recurrent sets play
an important role in the study of   qualitative behaviors of
dynamical systems. In a well-known paper   \cite{6}, Conley
discovered a fundamental connection between the chain recurrent
set and the collection of attractors for a deterministic dynamical
system $\varphi$  on a compact metric space $(X,d)$: the
complement of the chain recurrent set is the union of the sets
$B(A)-A$, where $A$ varies over the collection of attractors of
$\varphi$ and $B(A)$ is the basin of attraction of $A$. Namely,
$$X-CR(\varphi)=\bigcup [B(A)-A].$$
In other words, the complement of the chain recurrent set
$CR(\varphi)$ in the whole state space $X$ is the union of the
complements of attractors $A$ in their own basins of attraction.
This result is called  the Conley's theorem, which approaches the
fundamental theorem of dynamical systems \cite{8}. And the
fundamental theorem can be applied to such as the bifurcation
theorem\cite{Kap}.

The chain recurrent sets are widely studied and further extended by
many researchers in different contexts;  see for example
\cite{14,22,21,9,17,1,18}. In particular, by an alternative
definition of the chain recurrent set, Conley's theorem was
established for maps on locally compact metric spaces by
Hurley\cite{9,10} and was extended for semiflows and maps on
arbitrary metric spaces\cite{5,11}.
Let us recall a
definition from \cite{5,11}:\\

 Consider the set
$P(X)=\{f:X\rightarrow (0,\infty)\mid$ $f$ is continuous $\}$. For
any $\varepsilon \in P(X)$ and $T>0$, a finite sequence
$\{x_1,\cdots,x_n ,x_{n+1}; t_1, \cdots, t_n\}$, $x_i\in X$ and
$t_i\in  (-\infty,\infty), i=1,2,\cdots$, is called an
$(\varepsilon,T)$-chain for $\varphi$ from $p\in X$ to $q\in X$ if
$$x_1=p,\quad  x_{n+1}=q, \quad t_i \geq T,$$
and
$$d(\varphi(t_i,x_i),x_{i+1})<\varepsilon(\varphi(t_i,x_i)), \quad 1\leq i\leq n.$$

A point $p\in X$ is called chain recurrent for $\varphi$ if for
any $\varepsilon \in P(X)$ and $T>0$ there exists an
$(\varepsilon,T)$-chain for $\varphi$ from $p$ back to $p$. The
chain recurrent set $CR(\varphi)$ is the collection of all chain
recurrent points for $\varphi$.

A random dynamical system models dynamics in a state space $X$
influenced by probabilistic noise defined in a probability space
$(\Omega, \cF, \PX)$. Recently, the Conley's theorem was shown to
hold for a random dynamical system $\varphi$ on a compact metric
space $X$ by Liu \cite{1}. Liu obtained the following result by
defining a random chain recurrent set.

 \emph{Random Conley's theorem on compact metric space:}
 Let $\varphi$ be a random dynamical system on a  compact
space $X$. Assume that $U(\omega)$ is an arbitrary random
pre-attractor, $A(\omega)$ is the random local attractor
determined by $U(\omega)$, and $B(A)(\omega)$ is the  random basin
of attraction of $A(\omega)$. Then the following decomposition
holds:
$$X-CR_\varphi(\omega)=\bigcup [B(A)(\omega)-A(\omega)],\quad \mathbb{P}\mbox{-a.s}.$$
where the union is taken over all random local attractors
determined by random pre-attractors.  See Definition \ref{pre}
below for the meaning of pre-attractors and local attractors. Liu
\cite{2} further extended his   result to noncompact Polish spaces
(i.e., separable complete metric spaces) under an additional
\emph{absorbing condition}: the pre-attractor $U(\omega)$ is
assumed to be an absorbing set. Obviously this absorbing condition
depends on the underlying metric in the state space $X$.

In the present paper, we  prove a random Conley's theorem on
noncompact Polish metric spaces \emph{without the absorbing
condition}. Namely, we prove the following  result.

\begin{theorem} \label{main}
(Random Conley's theorem on noncompact Polish  spaces).\\
Let $\varphi$ be a random dynamical system on a non-compact Polish
space $X$. Assume that $U(\omega)$ is an arbitrary random
pre-attractor, $A(\omega)$ is the random local attractor
determined by $U(\omega)$ and $B(A,U)(\omega)$ is the basin of
attraction determined by $U(\omega)$ and $A(\omega)$. Then   the
following decomposition holds:
$$X-CR_\varphi(\omega)=\bigcup [B(A,U)(\omega)-A(\omega)],\quad \mathbb{P}\mbox{-a.s}.$$
where the union is taken over all random local attractors determined
by random pre-attractors.
\end{theorem}

Note that the \emph{random chain recurrent set} that we define
below (in Definition \ref{chain}) is   different from Liu's
\cite{1}, but they are equivalent when the state space $X$ is
compact.

This paper is organized as follows. After introducing
  random chain recurrent sets in \S 2,  we study the relation of random
chain recurrent sets and random attractors for random dynamical
systems in  \S 3, which includes the main result  of this paper.
Finally we demonstrate the main result by some simple examples in \S
4.

\section{Preliminaries} \label{prelim}

Throughout the    paper, we   assume that $(X,d)$  is a Polish
space, i.e. a separable complete metric space $X$ with metric $d$.
For $Y\subset X$, $Cl Y$, $int Y$, $Y^C$ denote respectively the
closure, interior and complement of $Y$.  Let $(\Omega,\mathcal
{F},\mathbb{P}) $ be a probability space, with sample space
$\Omega$, universal $\sigma$-algebra and probability measure
$\mathbb{P}$.

We recall the definition of a
 continuous random dynamical system (RDS) on the state space $X$, with
 the time set $\mathbb{T}= \mathbb{Z} $ or $\mathbb{R} $.
  For more information see references
\cite{4,19,1}.

 \begin{definition} (Random dynamical system).\\  A
 continuous random dynamical system (RDS), shortly denoted by
$\varphi$, consists of two ingredients: \\
(i) A model of the noise, namely a driving   flow $(\theta_t)_{t\in
\mathbb{T}}$ on the sample space $\Omega$, such that
$(t,\omega)\mapsto \theta_t \omega$ is a measurable flow that leaves
$\mathbb{P}$ invariant, i.e. $\theta_t \mathbb{P} =
\mathbb{P}$ for all $t \in \mathbb{T}$. \\
(ii) A model of the system influenced by noise, namely a cocycle
$\varphi$ over $\theta$, i.e. a measurable mapping $\varphi:
\mathbb{T} \times\Omega\times X\rightarrow  X, $
$(t,\omega,x)\mapsto \varphi(t,\omega,x)$, such that $(t,x)\mapsto
\varphi(t,\omega,x)$ is continuous for all $\omega \in \Omega$ and
the family $\varphi(t,\omega,\cdot) = \varphi(t,\omega) : X
\rightarrow X$ of random  mappings  satisfies the cocycle property:
\begin{eqnarray}
\varphi(0,\omega)=id_X , \varphi(t+s,\omega)=\varphi(t,\theta_s
\omega)\circ \varphi(s,\omega) \mbox{ for all } t,s \in \mathbb{T},
 \omega\in \Omega. \label{5}
\end{eqnarray}
 \end{definition}
It follows from (\ref{5}) that $\varphi(t,\omega)$ is a
homeomorphism of $X$ and
 $$\varphi(t,\omega)^{ -1 }=\varphi(-t,\theta_t\omega).$$

For a random variable $T(\omega)$, we call $T(\omega) > 0$ if it
holds $\mathbb{P}$-a.s. $\omega \in \Omega$.

\begin{definition} \label{pre}
(Pre-attractor and local attractor).\\
A random open set $U(\omega)$ is called a random pre-attractor of
a RDS $\phi$ if
\begin{eqnarray}
Cl{\bigcup \limits_{t\geq
\tau(\omega)}\varphi(t,\theta_{-t}\omega)
U(\theta_{-t}\omega)}\subset U(\omega)\quad \mbox{ for some }
\tau(\omega)>0,
\end{eqnarray}
where $\tau(\omega)$ is an $\mathcal {F}$-measurable random
variable.   We define the random local attractor $A(\omega)$
inside $U(\omega)$ to be the following:
\begin{eqnarray}
A(\omega)=\bigcap \limits_{n\in \mathbb{N}}Cl{\bigcup \limits_{t\geq
n\tau(\omega)}\varphi(t,\theta_{-t}\omega)
U(\theta_{-t}\omega)}.\label{7}
\end{eqnarray}
  The  basin of attraction for $A(\omega)$ is defined by
$B(A,U)(\omega)=\{x:\varphi(t,\omega)x\in U(\theta_t\omega)$ for
some $t\geq 0\}$.
\end{definition}

In the formula   \eqref{7}, we allow the random local attractor
$A(\omega)$  to be possibly an empty set. Similar to the  Lemma 3.1
in \cite{2}, we can assume that a pre-attractor $U(\omega)$ is a
forward invariant random open set. In the deterministic case, it can
be showed that if $X$ is a compact metric space and $U$ is a
nonempty open subset of $X$, then $U$ is a pre-attractor for a flow
$\varphi$ if and only if $Cl{\varphi(T,U)}\subset U$ for some $T>0$;
but when $X$ is not compact, this result may not be true \cite{11}.
When $X$ is compact, the basin of attraction of an attractor is
independent of the choice of the pre-attractor; but when $X$ is not
compact, the basin of attraction   generally depends on the
pre-attractor \cite{9}.

Now we define a random chain recurrent set in a Polish space. This
definition extends the definition in \cite{5,11,1}. Recall that a
function $f:X\times \Omega \rightarrow Y$ is said to be a
Caratheodory function if for all $x\in X$, $\omega \mapsto
f(x,\omega )$  is measurable, and for all $\omega \in \Omega$,
$x\mapsto f(x,\omega )$ is continuous. Let $Y=(0,+\infty)$ and we
use $ M(X\times \Omega)$ to denote the set of all such Caratheodory
functions.

\begin{definition} \label{chain}
(Random chain recurrent set).\\
(i) For a given Caratheodory function  $\varepsilon \in M(X\times
\Omega)$ and a random variable $T(\omega)>0$, a finite  sequence
$\{x_1(\omega), \cdots,x_n(\omega),x_{n+1}(\omega); t_1,\cdots,
t_n\} $ is called an $(\varepsilon,T)(\omega)$-chain for $\varphi$
from $x(\omega)$ to $y(\omega)$ if  for  $ 1\leq i\leq n$
  $$x_1(\omega)=x(\omega), \quad x_{n+1}(\omega)=y(\omega),\quad  t_i\geq T(\omega),$$
and
$$d(\varphi(t_i, \theta_{-t_i}\omega)x_i(\theta_{-t_i}\omega),x_{i+1}(\omega))
<\varepsilon(\varphi(t_i,
\theta_{-t_i}\omega)x_i(-\theta_{t_i}\omega),\omega).$$\\
(ii) Random variable $x(\omega)$ is called random chain recurrent
with index $\delta$ if
\begin{eqnarray*}
&&\mathbb{P}\{\ \omega\mid \mbox{ there exists an } (\varepsilon,
T)(\omega)\mbox{-chain beginning and ending at} \ x(\omega) \\&&
\mbox{ for any }\varepsilon \in M(X\times \Omega), \ T(\omega)> 0\}=
\delta,
\end{eqnarray*}
where index $\delta$ satisfies that for any $\eta> 0$, there exists
$\varepsilon_0 \in M(X\times \Omega)$ and $T_0(\omega)
> 0$  such that
any $(\varepsilon_0, T_0)(\omega)$-chain beginning and ending at $
x(\omega)$  with probability no more that $\delta+\eta$. Moreover, a
random variable $x(\omega)$ is called random chain recurrent if
$\delta=1$. A random variable  $x(\omega)$ is called partly random
chain recurrent if $1>\delta>0$.  And, $x(\omega)$ is
called completely random non-chain recurrent if  $\delta=0$. \\
\\
 (iii) $CR_\varphi(\omega)$ denotes the random chain recurrent
set of $ \varphi$, i.e.,
\begin{eqnarray*}
&&CR_\varphi(\omega)=\{x(\omega) \mid x(\omega) \mbox{ is the random
chain recurrent part of random }\\&&\mbox{ chain recurrent variable
with index}\ \delta\},
\end{eqnarray*}
where the random chain recurrent part is defined by
$x(\omega):\omega\in \Lambda$,
\begin{eqnarray*}
&&\Lambda:=\{\ \omega\mid \mbox{ there exists an } (\varepsilon,
T)(\omega)\mbox{-chain beginning and ending at} \ x(\omega)\\&&
\mbox{ for any }\varepsilon \in M(X\times \Omega), \ T(\omega)> 0\}.
\end{eqnarray*}
\end{definition}

The definition of completely  random non-chain recurrent variables
in the present paper is    different from that used in \cite{1}.

In the above definition, we   see that the random chain recurrent
set $CR_\varphi(\omega)$ of a RDS $\varphi$   has the property that
for any  random chain recurrent variable $x(\omega)$ with index
$\delta$, we have $x(\omega) \in CR_\varphi (\omega) $ with
probability $\delta$ and vice versa.

\begin{remark}
Notice that if $X$ is compact, then for any $\varepsilon\in
M(X\times\Omega)$, there exists $x_0\in X$ such that
$\varepsilon(x_0, \omega)=\min\limits_{x\in X}\{\varepsilon(x,
\omega)\}$. It can be shown that when $X$ is compact, the random
chain recurrent set that we define here is the same as   the one in
\cite{1}.
\end{remark}

\section{Random Conley's  theorem} \label{conley}

In this section we prove the main result of this paper, i.e.,
Theorem \ref{main}.  We first recall and prove some lemmas.

The following lemma about random local attractor is similar to
Lemma 3.1 in \cite{1}, which originates from \cite{19,20}.
\begin{lemma}\label{1}
Suppose $U(\omega)$ is a given pre-attractor and $ U(\tau(\omega)) =
Cl\bigcup\limits_{s\geq \tau(\omega)}$
$\varphi(s,\theta_{-s}\omega)U(\theta_{-s}\omega)$, then  the random
local attractor
\begin{eqnarray} A(\omega) = \bigcap \limits_{n\in \mathbb{N}}U(n \tau(\omega))
\end{eqnarray} determined by $U(\omega)$ is a  random closed set.
\end{lemma}

We prove a property of random chain recurrent variables.
\begin{lemma}\label{6}
If the random chain recurrent variable $x(\omega)\in U(\omega)$
$\mathbb{P}$-a.s. $\omega\in \Omega$, where $U(\omega)$ is a
random pre-attractor, then   $x(\omega)\in A(\omega)$
$\mathbb{P}$-a.s. $\omega\in \Omega$, where $A(\omega)$ is a
random attractor determined by $U(\omega)$.
\end{lemma}
{\bf Proof.} Recall that the pre-attractor $U(\omega) $ is forward
invariant. For $t\geq 0$, we have
\begin{eqnarray*}
\bigcup \limits_{t\geq 0}\varphi(t,\theta_{-t}\omega)
U(\theta_{-t}\omega)\subset U(\omega).
\end{eqnarray*}
We need to show that there exists an $\varepsilon \in M(X\times
\Omega)$ such that $\varepsilon \leq 1$ and
\begin{eqnarray}B(\varphi(t,\theta_{-t}\omega)x(\theta_{-t}\omega),
\varepsilon(\varphi(t,\theta_{-t}\omega)x(\theta_{-t}\omega),\omega))\subset
U(\omega)   \label{20}
\end{eqnarray}
 for all $x(\omega)\in U(\omega)$ and $t \geq 0$.

Let us construct such a Caratheodory function $\varepsilon$.
Define $\delta \in M(X\times \Omega)$ by
\begin{eqnarray*}
\delta(x,\omega)=\frac 1 2 \{d(x,{\bigcup \limits_{t\geq
0}\varphi(t,\theta_{-t}\omega) U(\theta_{-t}\omega}))+d(x,
X-U(\omega))\}.
\end{eqnarray*}
Then $\delta(x,\omega)>0$ since $x\notin X-U(\omega)$ if $x\in
{\bigcup \limits_{t\geq 0}\varphi(t,\theta_{-t}\omega)
U(\theta_{-t}\omega)}\subset U(\omega)$. Let $x(\omega)\in
U(\omega)$ and $t\geq 0$. For any $y(\omega)\in
B(\varphi(t,\theta_{-t}\omega)x(\theta_{-t}\omega),$
$\delta(\varphi(t,\theta_{-t}\omega)x(\theta_{-t}\omega),\omega))$,
\begin{eqnarray*}
d(\varphi(t,\theta_{-t}\omega)x(\theta_{-t}\omega),y(\omega))<\delta(\varphi(t,\theta_{-t}\omega)x(\theta_{-t}\omega),\omega)\\=\frac
1 2 d(\varphi(t,\theta_{-t}\omega)x(\theta_{-t}\omega),
X-U(\omega)).
\end{eqnarray*}
Thus we have
\begin{eqnarray*}
&&2d(\varphi(t,\theta_{-t}\omega)x(\theta_{-t}\omega),y(\omega))<d(\varphi(t,\theta_{-t}\omega)x(\theta_{-t}\omega),
X-U(\omega))\\& &\leq
d(\varphi(t,\theta_{-t}\omega)x(\theta_{-t}\omega),y(\omega))+d(y(\omega),
X-U(\omega)).
\end{eqnarray*}
Since $d(y(\omega),X-U(\omega))>
d(\varphi(t,\theta_{-t}\omega)x(\theta_{-t}\omega),y(\omega))\geq
0$, we have $y(\omega) \in U(\omega)$. Hence
\begin{eqnarray*}
B(\varphi(t,\theta_{-t}\omega)x(\theta_{-t}\omega),
\varepsilon(\varphi(t,\theta_{-t}\omega)x(\theta_{-t}\omega),\omega))\subset
U(\omega),
\end{eqnarray*}
where  $\varepsilon=\min\{\delta,1\}$ is the desired Caratheodory
function.

Let random variable $\widetilde{\tau}(\omega)> 0$ and $m$, $n$ be
positive integers. Since $x(\omega)\in CR_\varphi(\omega)$, there is
an $(\frac {\varepsilon} n , m\widetilde{\tau})(\omega)$-chain
$\{x_1(\omega),\cdots,x_k(\omega), x_{k+1}(\omega);t_1,$
$\cdots,$$t_k\}$ from $x(\omega)$ back to $x(\omega)$. Due to
\begin{eqnarray*}
d(\varphi(t_1,\theta_{-t_1}\omega
)x_1(\theta_{-t_1}\omega),x_2(\omega))<\frac 1 n
\varepsilon(\varphi(t_1,\theta_{-t_1}\omega)
x_1(\theta_{-t_1}\omega),\omega)\\
\leq \varepsilon(\varphi(t_1,\theta_{-t_1}\omega
)x_1(\theta_{-t_1}\omega),\omega),
\end{eqnarray*}
we see that
\begin{eqnarray*}
x_2(\omega)\in B(\varphi
(t_1,\theta_{-t_1}\omega)x(\theta_{-t_1}\omega),\varepsilon(\varphi
(t_1,\theta_{-t_1}\omega)x(\theta_{-t_1}\omega),\omega))\subset
U(\omega).
\end{eqnarray*}
Thus $x_k(\omega)\in U(\omega)$ by induction. Noticing that
\begin{eqnarray*}
d(\varphi(t_k,\theta_{-t_k}\omega
)x_k(\theta_{-t_k}\omega),x_{k+1}(\omega))<\frac 1 n
\varepsilon(\varphi(t_k,\theta_{-t_k}\omega)
x_k(\theta_{-t_k}\omega),\omega)\leq \frac 1 n,
\end{eqnarray*}
we obtain
\begin{eqnarray*}
d(x(\omega),\bigcup \limits_{t\geq
m\widetilde{\tau}(\omega)}\varphi(t,\theta_{-t}\omega)
U(\theta_{-t}\omega))\leq
d(x(\omega),\varphi(t_k,\theta_{-t_k}\omega
)x_k(\theta_{-t_k}\omega))< \frac 1 n.
\end{eqnarray*}
Let $n\rightarrow \infty$,  we conclude that
\begin{eqnarray*}
d(x(\omega),\bigcup \limits_{t\geq
m\widetilde{\tau}(\omega)}\varphi(t,\theta_{-t}\omega)
U(\theta_{-t}\omega) )=0.
\end{eqnarray*}
 Hence $x(\omega)\in Cl{\bigcup \limits_{t\geq
m\widetilde{T}(\omega)}\varphi(t,\theta_{-t}\omega)
U(\theta_{-t}\omega)} $. This implies that
\begin{eqnarray*}
x(\omega)\in \bigcap \limits_{m\in \mathbb{N}} Cl{\bigcup
\limits_{t\geq
m\widetilde{\tau}(\omega)}\varphi(t,\theta_{-t}\omega)
U(\theta_{-t}\omega)}=A(\omega).
\end{eqnarray*}
 \hfill$\square$

Finally we   obtain another property for random chain recurrent
variables.
\begin{lemma}\label{10}
If the random chain recurrent variable $x(\omega)\in B(A,
U)(\omega)$ $\mathbb{P}$-a.s. $\omega\in \Omega$, then
$x(\omega)\in A(\omega)$ $\mathbb{P}$-a.s. $\omega\in \Omega$,
where $B(A,U)(\omega)$ is the basin of attraction determined by
$U(\omega)$ and $A(\omega)$.
\end{lemma}
{\bf Proof.}  
Let a sequence
$s_n\rightarrow \infty$, denote
\begin{eqnarray*}
&&\Omega_n=\{\omega\mid \varphi(s_n,\omega)x(\omega)\in
U(\theta_{s_n}\omega)\},\\
&&\Omega'_n=\{\omega\mid
\varphi(s_n,\theta_{-s_n}\omega)x(\theta_{-s_n}\omega)\in
U(\omega)\}.
\end{eqnarray*}
Then we have
$$\mathbb{P}(\Omega)=\mathbb{P}(\bigcup\limits_{n=1}^\infty\Omega_n)=
\lim\limits_{n\rightarrow
\infty}\mathbb{P}(\Omega_n)=\lim\limits_{n\rightarrow
\infty}\mathbb{P}(\Omega_n').$$ Since $\mathbb{P}(\Omega_n')\leq
\mathbb{P}(\bigcup\limits_{i=1}^n\Omega_i')$, we have
$\mathbb{P}(\bigcup\limits_{n=1}^\infty\Omega_n')=1$. For fixed
$\omega \in \bigcup\limits_{i=1}^n\Omega_i'$ and $x(\omega)\in
CR_\varphi(\omega)$, there exists $s_n$ such that the following
holds:
\begin{eqnarray*}
\varphi(s_n,\theta_{-s_n}\omega)x(\theta_{-s_n}\omega) \in
U(\omega).
\end{eqnarray*}

Assume that random variable $\widetilde{\tau}(\omega)>0$ and $m$,
$n$ are positive integers. Select $\varepsilon \in M(X\times
\Omega)$ satisfies \eqref{20} for
$\varphi({s_n},\theta_{-{s_n}}\omega)x(\theta_{-{s_n}}\omega)\in
U(\omega)$. Since $x(\omega)\in CR_\varphi(\omega)$, there is an
$(\frac {\varepsilon} n , m\widetilde{\tau})(\omega)$-chain
$\{x_1(\omega),\cdots,x_k(\omega),$ $ x_{k+1}(\omega);$$ t_1,$ $
\cdots,$ $t_k\}$ from $x(\omega)$ back to $x(\omega)$. Take
$t_1>{s_n}$ by $m$ large enough. Since
\[ \begin{split}
&d(\varphi(t_1-{s_n},\theta_{-{(t_1-{s_n})}}\omega )\circ
\varphi({s_n},\theta_{-{s_n}}\theta_{{s_n}-t_1}\omega)x_1(\theta_{-{s_n}}\theta_{{s_n}-t_1}\omega),
x_2(\omega))\\
&=d(\varphi(t_1,\theta_{-{t_1}}\omega
)x_1(\theta_{-{t_1}}\omega),x_2(\omega))\\
&<\frac 1 n \varepsilon(\varphi(t_1,\theta_{-{{{t_1}}}}\omega)
x_1(\theta_{-{t_1}}\omega),\omega)\\
&\leq  \varepsilon(\varphi(t_1,\theta_{{{-{t_1}}}}\omega
)x_1(\theta_{-{t_1}}\omega),\omega).
\end{split} \]
We have
\begin{eqnarray*}
x_2(\omega)\in  B(\varphi
({t_1},\theta_{-{t_1}}\omega)x(\theta_{-{t_1}}\omega),\varepsilon(\varphi
(t_1,\theta_{-{t_1}}\omega)x(\theta_{-{t_1}}\omega),\omega)) \subset
U(\omega).
\end{eqnarray*}
By a similar argument to the proof of Lemma \ref{6}. We have
\begin{eqnarray*}
x(\omega)\in \bigcap \limits_{m\in \mathbb{N}} Cl{\bigcup
\limits_{t\geq
m\widetilde{\tau}(\omega)}\varphi(t,\theta_{-t}\omega)
U(\theta_{-t}\omega)}=A(\omega).
\end{eqnarray*}
 \hfill$\square$

 \begin{remark}  The   proof of Lemma \ref{10} is
 different from that   in \cite{1}. In \cite{1}, the  entrance time
 of $x(\omega)$ into $U(\omega)$ under the cocycle $\varphi$ does not depend on the $\omega$
 finally, while the entrance time depends on the $\omega$ in this
Lemma.
\end{remark}

 \begin{remark} \label{8}
For a given random variable $x(\omega)$, define
$\Omega_{CR}(x)=\{\omega\mid x(\omega)\in CR_\varphi(\omega)\}$.
 If $x(\omega)$ is partly random chain recurrent and $x(\omega)\in B(A,U)(\omega)$
$\mathbb{P}$-a.s. $\omega\in \Omega$,  we know that if $\omega \in
\Omega_{CR}(x)$, then
\begin{eqnarray}
x(\omega) \in A(\omega) \label{12}
\end{eqnarray}
holds  for any given random local attractor $A(\omega)$. This can
be seen as in the proof   of Lemma \ref{10}. Hence by a similar
discussion of \cite{1} we   obtain that
\begin{eqnarray*}
\bigcup [B(A,U)(\omega)-A(\omega)] \subset X - CR_\varphi (\omega),
\quad\mathbb{P}\mbox{-a.s.}
\end{eqnarray*}
\end{remark}
\medskip

We are now ready to prove our main result.\\ \\

\noindent {\bf Proof of Theorem \ref{main}:} \\

 Suppose $x(\omega)$ is a
random variable, $\varepsilon\in M(X\times \Omega)$, random variable
$\tau(\omega)>0$. Define
 \begin{eqnarray*}
 U_1(\omega):&=&\bigcup \limits_{t\geq
 \tau(\omega)}B(\varphi(t,\theta_{-t}\omega)x(\theta_{-t}\omega),\varepsilon(\varphi(t,\theta_{-t}\omega)x(\theta_{-t}\omega),\omega)),\\
 U_2(\omega):&=&\bigcup \limits_{t\geq \tau(\omega)}\bigcup \limits_{y(\omega)\in
 U_1(\omega)}B(\varphi(t,\theta_{-t}\omega)y(\theta_{-t}\omega),\varepsilon(\varphi(t,\theta_{-t}\omega)y(\theta_{-t}\omega),\omega)),\\
             &\cdots&\\
 U_n(\omega):&=&\bigcup \limits_{t\geq \tau(\omega)}\bigcup \limits_{y(\omega)\in
 U_{n-1}{(\omega)}}B(\varphi(t,\theta_{-t}\omega)y(\theta_{-t}\omega),\varepsilon(\varphi(t,\theta_{-t}\omega)y(\theta_{-t}\omega),\omega)),\\
             &\cdots.&
 \end{eqnarray*}
 For fixed $x\in X$, define
\begin{eqnarray*}
d_1(t,\omega):=d(x,B(\varphi(t,\theta_{-t}\omega)x(\theta_{-t}\omega),\varepsilon(\varphi(t,\theta_{-t}\omega)x(\theta_{-t}\omega),\omega)))
\end{eqnarray*}
Since $\varphi(t,\theta_{-t}\omega)x(\theta_{-t}\omega)$,
$\varepsilon(\varphi(t,\theta_{-t}\omega)x(\theta_{-t}\omega),\omega)$
is $\mathcal {B}(\mathbb{T})\times \mathcal {F}$-measurable,  we can
obtain $(t,\omega)\rightarrow d_1(t,\omega)$ is $\mathcal
{B}(\mathbb{T})\times \mathcal {F}$-measurable.
\begin{eqnarray*}
d(x,U_1(\omega))&=&d(x,\bigcup \limits_{t\geq
\tau(\omega)}B(\varphi(t,\theta_{-t}\omega)x(\theta_{-t}\omega),\varepsilon(\varphi(t,\theta_{-t}\omega)x(\theta_{-t}\omega),\omega)))\\&=&\inf\limits_{t\geq
\tau(\omega)} d_1(t,\omega).
\end{eqnarray*}
For arbitrary $a\in R^+$, we have
\begin{eqnarray*}
\{\omega\mid \inf\limits_{t\geq
\tau(\omega)}d_1(t,\omega)<a\}=\Pi_\Omega \{(t,\omega)\mid
d_1(t,\omega)<a, t\geq \tau(\omega)\}.
\end{eqnarray*}
It is obvious that the function $(t,\omega)\rightarrow
t-\tau(\omega)$ is measurable with respect to $\mathcal
{B}(\mathbb{T})\times \mathcal {F}$, so by the Projection Theorem,
we obtain that $U_1(\omega)$ is a random open set. Similarly we can
prove that $U_2(\omega),\cdots,$ $ U_n(\omega),$ $ \cdots$ are all
 random open sets. So the set
\begin{eqnarray}
U_x(\omega):=\bigcup \limits_{n\in \mathbb{N}}U_n(\omega)\label{4}
\end{eqnarray} is a random open set.
From the construction of $U_x(\omega)$ we can see that $U_x(\omega)$
is the set of all possible end points of $(\varepsilon,
\tau)(\omega)$-chains that begin at $x(\omega)$.

In the following we prove that  $U_x(\omega)$ is a pre-attractor. It
is easy to see that $x(\omega)\notin U_x(\omega)$ and
 \begin{eqnarray*}
 \varphi(t,\theta_{-t}\omega)x(\theta_{-t}\omega)\in U_x(\omega)
 \end{eqnarray*}
 when $t\geq \tau(\omega)$.
Since $\varepsilon(x,\omega)$ is continuous in $x$, there exists a
map $\delta: X\times \Omega\rightarrow (0,+\infty)$ with $\delta\leq
\frac{\varepsilon} 2$ such that $\varepsilon(y(\omega),\omega)>\frac
1 2 \varepsilon(x(\omega),\omega)$ when
$d(x(\omega),y(\omega))<\delta(x(\omega),\omega)$. For
$B(y(\omega),\delta(y(\omega),\omega))\bigcap
\varphi(t,\theta_{-t}\omega)U(\theta_{-t}\omega)\neq \emptyset$,
$t\geq \tau(\omega)$. There exists $z(\omega)\in U_x(\omega)$ such
that
\begin{equation*}
  d(\varphi(t,\theta_{-t}\omega)z(\theta_{-t}\omega),y(\omega))<\delta(y(\omega),\omega),
  \quad   t\geq \tau(\omega).
\end{equation*}
Since
$d(\varphi(t,\theta_{-t}\omega)z(\theta_{-t}\omega),y(\omega))<\delta(y(\omega),\omega)$,
we have
$$\varepsilon(\varphi(t,\theta_{-t}\omega)z(\theta_{-t}\omega),\omega)>\frac
1 2 \varepsilon(y(\omega),\omega).$$ Thus
 {\setlength\arraycolsep{2pt}
\begin{eqnarray*}
  d(\varphi(t,\theta_{-t}\omega)z(\theta_{-t}\omega),y(\omega))&<&\delta(y(\omega),\omega)\leq \frac 1 2
  \varepsilon(y(\omega),\omega)\\&<&
  \varepsilon(\varphi(t,\theta_{-t}\omega)z(\theta_{-t}\omega),\omega).
\end{eqnarray*}}
Hence there exists an $(\varepsilon,\tau)(\omega)$-chain from
$x(\omega)$ to $y(\omega)$. This means that
$$Cl{\bigcup\limits_{t\geq
\tau(\omega)}\varphi(t,\theta_{-t}\omega)U(\theta_{-t}\omega)}\subset
U(\omega).$$
 So $U_x(\omega)$ is a random pre-attractor and it
determines a random closed local attractor $A_x(\omega)$ by Lemma
\ref{1}.

 When $x(\omega)$ is  random chain recurrent
with index $\delta$, i.e. $x(\omega) \in  X - CR_\varphi(\omega) $
with probability $1-\delta$, by a similar argument   in \cite{1},
we obtain
$$ X - CR_\varphi(\omega) \subset \bigcup
[B(A,U)(\omega)-A(\omega)], \quad \mathbb{P}\mbox{-a.s.} $$

By Remark \ref{8}, we complete the proof of Theorem \ref{main}.
\hfill$\square$\\


\section{Examples}

We now present some   examples to illustrate our result in Theorem
\ref{main}. It is known that when the pre-attractor is bounded or
the state space is finite dimensional, then the pre-attractor is an
absorbing set. The pre-attractor in Example \ref{2.3} is unbounded
and it is not an absorbing set. The pre-attractor in Example
\ref{2.1} is in infinite dimensional state space and it is not
necessarily an absorbing set. In both examples, our Theorem
\ref{main} applies, while Liu's result \cite{2} does not apply as it
requires  the pre-attractor to be also an absorbing set.


\begin{example}\label{2.3}
Consider a random map $f:((x,y),\omega)\rightarrow
((x+1,y),\omega)$.   It generates a discrete random dynamical
system on the state space $X=\mathbb{R}^2$. The set
  $$U(\omega)=\{(x,y)\mid |y|\leq e^x\}$$ is a pre-attractor. But it
  does not satisfy the absorbing condition. It is clear that the local attractor is $A(\omega)=\{(x,0)\mid x\in
  \mathbb{R}\}$, while $B(A,U)(\omega)=\mathbb{R}^2$ is the basin of attraction.  By Theorem \ref{main},
  we   conclude that
  the random chain recurrent set is $\emptyset$.
\end{example}

We revise an example from \cite{Crau} to fit our purpose here.
\begin{example}\label{2.1}
 Let $X$ be a Polish space and denoted  $Pr(X)$  all Borel
 probability measures on $X$. Suppose that $\varphi$ is a  $RDS$ on
  $X$. Then $\varphi$ induces a  RDS on the Polish
 space $Pr(X)$, as follows
 \begin{eqnarray*}
   \Phi: \mathbb{T}\times\Omega \times Pr(X)&\rightarrow &Pr(X)\\
(t,\omega,\rho)&\mapsto& \varphi(t,\omega)\rho,
 \end{eqnarray*}
 where $\varphi(t,\omega)\rho$ is the image measure of the
 deterministic measure $\rho$ under the random map
 $\varphi(t,\omega)$.

 Assume that there exists a unique pre-attractor
  $\Gamma_U\neq
  Pr(X)$ and  $Pr(A(\omega))$ is the corresponding local
  attractor. Then by Theorem \ref{main},
   $X-CR_\Phi(\omega)=B(Pr(A),\Gamma_U)(\omega)-Pr(A(\omega)),\quad \mathbb{P}\mbox{-a.s.}$

\end{example}

\end{document}